\documentclass[11pt]{amsart}
\usepackage[T2A,T1]{fontenc}
\usepackage[utf8]{inputenc}		
\usepackage{ucs}
\usepackage[russian,english]{babel}	
\usepackage{amssymb}
\usepackage{ dsfont }
\usepackage{multirow}
\usepackage{amsthm}
\usepackage{tikz}
\usepackage{amsmath}
\usepackage{graphicx}
\usepackage{wrapfig}
\usepackage{caption}
\usepackage[left=2.5cm,right=2.5cm,top=2.4cm,bottom=2.4cm,bindingoffset=0cm]{geometry}
\usepackage{enumerate}
\usepackage{makecell}	
\newtheorem{state}{Statement}[section]
\newtheorem{theorem}{Theorem}
\newtheorem{conj}{Conjecture}

\newtheorem{cor}{Corollary}

\newtheorem{prob}{Problem}

\newtheorem{prop}[state]{Proposition}

\makeatletter
\newcommand*{\rom}[1]{\expandafter\@slowromancap\romannumeral #1@}
\makeatother

\newcommand{\Aut}{\mathop{\mathrm{Aut}}\nolimits}

\makeatletter
\def\@seccntformat#1{\csname the#1\endcsname. }
\def\@biblabel#1{#1.}

\title[On arithmetical properties of finite groups]{On arithmetical properties and arithmetical characterizations of finite groups}
\author{Natalia~V.~Maslova}
\address{Natalia Vladimirovna Maslova\newline Krasovskii Institute of Mathematics and Mechanics UB RAS,\newline
16, S. Kovalevskaja str.,
Yekaterinburg, 620108, Russia\newline
Ural Federal University,\newline
19, Mira str., Yekaterinburg, 620002, Russia \newline
ORCID: 0000-0001-6574-5335}
\email{butterson@mail.ru}

\begin{document}


\maketitle

\hfill To Professor Viktor Danilovich Mazurov on the occasion of his 80th birthday

\hfill and to Professor Wujie Shi on the occasion of his 80th birthday

\medskip

This paper is based on the talk given by the author at the conference $(WM)^2$~--- World Meeting for Women in Mathematics~---$2022$.

\medskip


\medskip

MSC classes: 20-02, 20D60, 20D05, 20D06, 20D08, 05C15, 05C25, 05C99, 05E30

\section{Introduction}

Symmetry is one of the fundamental principles of self-organization in nature. By studying the group of symmetries of an object (for example, a graph, a crystal, a dynamical  system, a network, and so on), we can get new information about the object itself.
However the situation when the group of symmetries of an object is known {\it a priori}, is rather rare. But from some 'visible' properties of the object, it can be possible to extract, for example, information about arithmetical properties of its group of symmetries (i.e., the properties of the group which are defined by its arithmetical parameters such as the order of the group, the element orders and so on). The problem of defining a group or describing at least some of its structural properties and features of possible actions on objects, if only some arithmetic parameters of this group are known, naturally appears (see Figure~1).

\begin{center}
\includegraphics[width=0.65\textwidth]{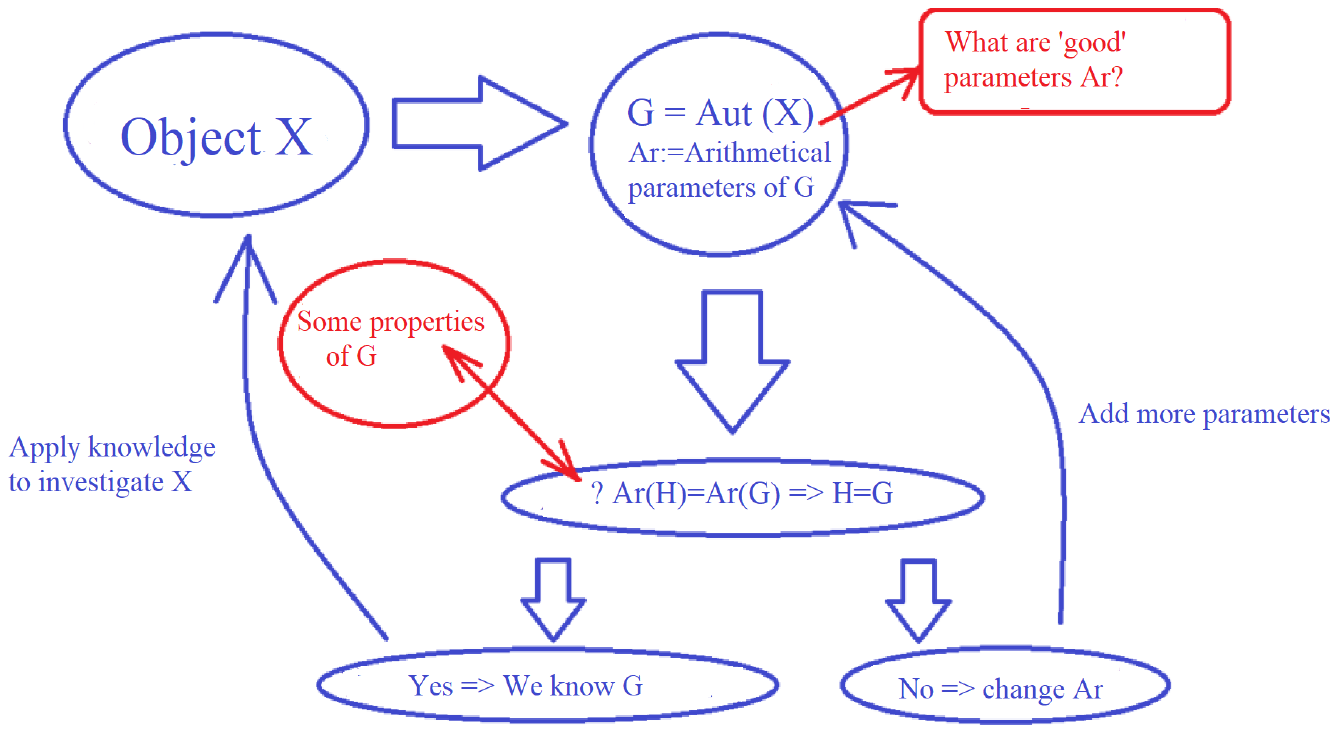}\\
Figure~1
\end{center}

The problem of characterization of a (finite) group by its arithmetical parameters is actively investigated in modern group theory, and obtaining results of this kind is a development of a mathematical apparatus that can later be applied outside of mathematics. In this paper, we discuss a number of results on arithmetical properties of finite groups and characterizations of finite groups by their arithmetical parameters obtained recently.

\medskip

Throughout the paper we consider only finite groups and simple graphs, and henceforth the term {\it group} means finite group, the term {\it graph} means simple graph (undirected graph without loops and multiple edges).

\smallskip

Let $\Omega$ be a finite set of positive integers. Define $\pi(\Omega)$ to be the set of all prime divisors of integers from $\Omega$.
If $\Omega=\{n\}$, then we write $\pi(n)$ instead of $\pi(\{n\})$. A graph $\Gamma(\Omega)$ whose vertex set is $\pi(\Omega)$ and two distinct vertices~$p$ and~$q$ are adjacent if and only if $pq$ divides some element from $\Omega$ is called the {\it prime graph} of $\Omega$.

\medskip

{\bf Example.} If $\Omega=\{2,15,21,35\}$, then $\pi(\Omega)=\{2,3,5,7\}$ and $\Gamma(\Omega)$ is as in Figure~2.

\medskip

\begin{center}    \begin{tikzpicture}
        \tikzstyle{every node}=[draw,circle,fill=white,minimum size=4pt, inner sep=0pt]
        \draw (0,0) node (2) [label=below:$2$] {}
        (0:1.0cm) node (3) [label=below:$3$] {}
        (0:2.0cm) node (5) [label=right:$5$] {}
        (-20:2.0cm) node (7) [label=below:$7$] {}
        (3)--(5)
        (3)--(7)
        (5)--(7);
    \end{tikzpicture}
    \\
Figure~2
    \end{center}

\smallskip

Let~$G$ be a group. Denote by  $\omega (G)$ the {\it spectrum} of~$G$, that is, the set of all its element orders. Let $\pi(G)=\pi(|G|)$ be the {\it prime spectrum} of $G$. By Lagrange's theorem and Cauchy's theorem, $\pi(G)=\pi(\omega(G))$. The set $\omega(G)$ defines the {\it Gruenberg--Kegel graph} (or the {\it prime graph}; or {\it prime vertex graph} as in~\cite{Iranmanesh_Praeger}; or {\it GK graph} as in~\cite{Cameron}) $\Gamma(G)$ of~$G$ by $\Gamma(G)=\Gamma(\omega(G))$, i.\,e., in this graph the vertex set is $\pi(G)$, and distinct vertices~$p$ and~$q$ are adjacent if and only if $pq\in\omega (G)$. In many papers the authors use the term {\it prime graph} for various reasons, mainly for brevity. But in my opinion it seems more correct to use the term {\it Grunberg-Kegel graph} to distinguish this graph from other well-known kinds of graphs such as {\it character degree graphs} (see, for example,~\cite{Akhlaghi_etal}) which are also defined as prime graphs built with other sets of arithmetical parameters of groups.

\medskip

The concept of Gruenberg--Kegel graph appeared in the unpublished manuscript~\cite{Gruen_Keg} by K.~Gruenberg and O.~Kegel and proved to be very useful, for example, with connection to research on some cohomological questions in integral group rings: the augmentation ideal of an integral group ring is decomposable as a module if and only if the Gruenberg--Kegel graph of the group is disconnected~\cite{Gruen_Rogg}. Also in~\cite{Gruen_Keg}, the groups with disconnected Gruenberg--Kegel graph where characterized; this result was published later in the paper~\cite{Williams} by J.~Williams, who was a student of K.~Gruenberg, and now this theorem is well-known as the Gruenberg--Kegel Theorem (see Theorem~\ref{Gruenberg--Kegel theorem} below). Connections of the Gruenberg--Kegel graph with some graphs defined on groups were discussed in~\cite{Cameron}.

\medskip

It is easy to see that for groups $G$ and $H$, if $G \cong H$, then $\omega(G)=\omega(H)$; and if $\omega(G)=\omega(H)$, then $\Gamma(G)=\Gamma(H)$; if $\Gamma(G)=\Gamma(H)$, then $\pi(G)=\pi(H)$ and $\Gamma(G)$ and $\Gamma(H)$ are isomorphic as abstract graphs (i.\,e., as unlabelled graphs). The converse does not hold in each case, as the following series of examples demonstrates (see, for example, \cite{Atlas}):

\begin{itemize}

\item $S_5\not \cong S_6$ but $\omega(S_5)=\omega(S_6)${\rm;}

\item $\omega(A_5)\not =\omega(A_6)$ but $\Gamma(A_5)=\Gamma(A_6)${\rm;}

\item $\Gamma(A_{10})\not =\Gamma(\Aut(J_2))$ but  $\Gamma(A_{10})$ and $\Gamma(\Aut(J_2))$ are isomorphic as abstract graphs  and $\pi(A_{10}) =\pi(\Aut(J_2))$, see Figure~3.
    \begin{center}
    \begin{tikzpicture}
    \tikzstyle{every node}=[draw,circle,fill=white,minimum size=4pt,
                            inner sep=0pt]
    \draw (0,0) node (2) [label=left:$2$]{}
        -- ++ (-30:1cm) node (3) [label=below:$3$]{}
        -- ++ (0:1cm) node (7) [label=below:$7$]{}
          (-90:1.0cm) node (5) [label=left:$5$]{};
    \draw (5) -- (3);
    \draw (5) -- (2);
    \end{tikzpicture}{$\Gamma(A_{10})$;} $\mbox{ }$$\mbox{ }$$\mbox{ }$$\mbox{ }$$\mbox{ }$$\mbox{ }$$\mbox{ }$$\mbox{ }$
    \begin{tikzpicture}
    \tikzstyle{every node}=[draw,circle,fill=white,minimum size=4pt,
                            inner sep=0pt]
    \draw (0,0) node (2) [label=left:$3$]{}
        -- ++ (-30:1cm) node (3) [label=below:$2$]{}
        -- ++ (0:1cm) node (7) [label=below:$7$]{}
          (-90:1.0cm) node (5) [label=left:$5$]{};
    \draw (5) -- (3);
    \draw (5) -- (2);
    \end{tikzpicture}{$\Gamma(\Aut(J_2))$.}\\
    Figure~3
\end{center}

\bigskip

\end{itemize}

\noindent A group $G$ is
\begin{itemize}
\item {\it recognizable} by spectrum (respectively, by Gruenberg--Kegel graph) if for each group~$H$, $\omega(G)=\omega(H)$ (respectively, $\Gamma(G)=\Gamma(H)$) if and only if $G \cong H${\rm;}
\item  {\it $k$-recognizable} by spectrum (respectively, by Gruenberg--Kegel graph), where $k$ is a positive integer, if there are exactly $k$ pairwise non-isomorphic groups  $H$ with $\omega(H)=\omega(G)$ (respectively, $\Gamma(H)=\Gamma(G)$){\rm;}
\item  {\it almost recognizable} by spectrum (respectively, by Gruenberg--Kegel graph) if it is $k$-recognizable by spectrum (respectively, by Gruenberg--Kegel graph) for some positive integer $k${\rm;}
\item {\it unrecognizable} by spectrum (respectively, by Gruenberg--Kegel graph), if there are infinitely many pairwise non-isomorphic groups $H$ with $\omega(H)=\omega(G)$ (respectively, $\Gamma(H)=\Gamma(G)$).
\end{itemize}

Note that if a group $G$ is recognizable by Gruenberg--Kegel graph, then $G$ is recognizable by spectrum. The converse does not hold since, for example, $\Gamma(A_5)=\Gamma(A_6)$, but the group $A_5$ is recognizable by spectrum (see~\cite{Shi86}) while the group $A_6$ is unrecognizable by spectrum (\cite{Brandl_Shi}, see also \cite[Theorem~2]{Lytkin}), therefore both groups $A_5$ and $A_6$ are unrecognizable by Gruenberg--Kegel graph.

\medskip

Recall that a group $G$ is {\it simple} if $G$ does not contain proper normal subgroups and is {\it almost simple} if there exists a non-abelian simple group $S$ such that $$S \cong Inn(S) \unlhd G \le Aut(S);$$ in this case, $Soc(G)\cong S$ (the {\it socle} $Soc(G)$ is the subgroup of $G$ generated by all its minimal non-trivial normal subgroups).  In accordance with the Classification of Finite Simple Groups (CFSG; the reader can find out more about this in~\cite{Gorenstein}), non-abelian simple groups are contained in the following list:

\begin{itemize}

\item {alternating groups}{\rm;}

\item {classical groups}{\rm:} $PSL_n(q)$, $PSU_n(q)$, $PSp_{2n}(q)$, $P\Omega_{2n+1}(q)$, $P\Omega_{2n}^+(q)$, $P\Omega_{2n}^-(q)${\rm;}

\item {exceptional groups of Lie type}{\rm:} ${^2}B_2(2^{2n+1})$, $G_2(q)$, ${^2}G_2(3^{2n+1})$, ${^3}D_4(q)$, $F_4(q)$, ${^2}F_4(2^{2n+1})$, $E_6(q)$, ${^2}E_6(q)$, $E_7(q)$, $E_8(q)${\rm;}

\item $26$ {sporadic} groups.

\end{itemize}

\medskip

A non-abelian simple group $L$ is

\begin{itemize}

\item {\it quasirecognizable} by spectrum (respectively, by Gruenberg--Kegel graph), if every group $H$ with $\omega(G)=\omega(H)$ (respectively, $\Gamma(G)=\Gamma(H)$) has a unique nonablelian composition factor $S$ and $S\cong L${\rm;}

\item {\it recognizable by spectrum among covers} (respectively, {\it recognizable by Gruenberg--Kegel graph among covers}) if for each group $H$, $H/N\cong L$ and $\omega(H)=\omega(L)$  (respectively, $H/N\cong L$ and $\Gamma(H)=\Gamma(L)$) if and only if $H\cong L${\rm;}

\item {\it recognizable by spectrum among automorphic extensions} (respectively, {\it recognizable by Gruenberg--Kegel graph among automorphic extensions}) if for each group $H$, $Soc(H)=L$ and $\omega(H)=\omega(L)$ (respectively, $Soc(H)=L$ and $\Gamma(H)=\Gamma(L)$) if and only if $H\cong L$.

\end{itemize}

\medskip

\medskip

In this paper we will discuss some results on characterization of a group by spectrum and by Gruenberg--Kegel graph and also some results on combinatorial properties of Gruenberg--Kegel graphs. Note that recently W.~Shi published a survey paper~\cite{Shi_Surv2024} in which he discusses characterizations of finite groups (mostly simple groups) by some other arithmetical parameters.

\section{Characterization by spectrum}

If $G$ is a group then denote by $S(G)$ the {\it solvable radical} of $G$, i.\,e., the largest solvable normal subgroup of $G$.

\smallskip
The following criterion of recognition of a group by spectrum was obtained by V.~D.~Mazurov and W.~Shi~\cite{Mazurov_Shi}.

\begin{theorem}[{\cite[Theorem]{Mazurov_Shi}}]\label{UnrecSpec} Let $G$ be a group. The following statements are equivalent{\rm:}
\begin{itemize}
\item[$(1)$] $G$ is unrecognizable by spectrum{\rm;}
\item[$(2)$] there exists a group $H$ such that $S(H)\not =1$ and $\omega(H)=\omega(G)$.
\end{itemize}
\end{theorem}

\begin{cor} If $G$ is almost recognizable {\rm(}in particular, recognizable{\rm)} by spectrum, then $S(G)=1$. Thus, the question of recognizability by spectrum is of interest for almost simple groups in particular for non-abelian simple groups.
\end{cor}

The first examples of groups which are recognizable by spectrum were found in the middle of the 1980s by W.~Shi~\cite{Shi84,Shi86}. Since then the problem of recognition by spectrum was very actively investigated. A survey of this research area can be found in the paper by M.~A.~Grechkoseeva, V.~D.~Mazurov, W.~Shi, A.~V.~Vasil'ev, and N.~Yang~\cite{GreMazShiVasYang} the first four of whom made very significant contributions to investigations of characterization of finite groups by spectrum. Also we refer the reader to earlier surveys by V.~D.~Mazurov~\cite{Mazurov2005,Mazurov2015}. In particular, the following theorem is true.

\begin{theorem}[{\cite[Theorem~2.1]{GreMazShiVasYang}}]\label{SpectrumMain} Let $L$ be one of the following non-abelian simple groups:

\begin{itemize}

\item a sporadic group other than $J_2$;

\item an alternating group $A_n$, where $n \not \in \{6, 10\}$;

\item an exceptional group of Lie type other than ${^3}D_4(2)$;

\item classical group of a sufficiently large dimension.

\end{itemize}
Then every group $H$ such that $\omega(H)=\omega(L)$ is isomorphic to some group $G$ with $L \unlhd G \le Aut(L)$.
In particular, $L$ is almost recognizable by spectrum.

\end{theorem}

Moreover, in~\cite{Shi89}, W.~Shi made the following conjecture (the history of this conjecture is discussed in \cite[p.~170]{GreMazShiVasYang}).

\begin{conj}[{\cite[Conjecture]{Shi89}}]\label{OrdSpecConj} If $G$ is a simple group and $H$ is a group, then $H\cong G$ if and only if $|H|=|G|$ and $\omega(H)=\omega(G)$.
\end{conj}

In~\cite{Shi89}, Conjecture~\ref{OrdSpecConj} was proved for all sporadic simple groups. Later A.~S.~Kondrat'ev added Conjecture~\ref{OrdSpecConj} to the Kourovka Notebook~\cite[Problem~12.39]{Kourovka}. The conjecture was investigated in a series of papers by W.~Shi, J.~Bi, H.~Cao, and M.~Xu~\cite{Cao_Shi_Order5,Shi_Order3,Shi_Bi_Order4,Shi_Bi_Order2,Shi_Bi_Order1,Xu_Shi_Order6}, and finally the proof of Conjecture~\ref{OrdSpecConj} was completed in 2009 by A.~V.~Vasil'ev, M.~A.~Grechkoseeva, and V.~D.~Mazurov~\cite{VasGrechMaz09}. Thus, the following theorem is true.

\begin{theorem}[{\cite[Theorem]{VasGrechMaz09}}]\label{OrdSpecThm} If $G$ is a simple group and $H$ is a group such that $|H|=|G|$ and $\omega(H)=\omega(G)$, then $H \cong G$.
\end{theorem}

Some almost simple but not simple groups are known to be recognizable by spectrum. For example, if $n\not =10$, then the recognition question by spectrum has been solved for the symmetric group $S_n$, and if $n$ is sufficiently large, then $S_n$ is recognizable by spectrum~\cite{Gorsh_Grish}. Note that one of the first results in this direction was obtained by C.~E.~Praeger and W.~Shi~\cite{Praeger_Shi}, they proved that the group $S_7$ is recognizable by spectrum but the group $S_8$ is unrecognizable by spectrum.

\smallskip

Moreover, there are groups with non-simple socle which are recognizable by spectrum. The first such examples were found by V.~D.~Mazurov~\cite{Mazurov97}, who proved that the group ${^2}B_2(2^7)\times {^2}B_2(2^7)$ and the permutation wreath product of ${^2}B_2(2^7)$ with $H \cong 23:11 < S_{23}$ are recognizable by spectrum. The next result of this kind was obtained only after $24$ years by I.~B.~Gorshkov and the author in~\cite{Gorsh_Mas}, where we proved that the direct product $J_4 \times J_4$ of two copies of the simple sporadic group $J_4$ is recognizable by spectrum.

Note that if the direct product of $k$ copies of a group $G$ is recognizable by spectrum, then for each $i \le k$ the direct product of $i$ copies of $G$ is recognizable by spectrum. Thus, the following problems are of interest.

\begin{prob}[{\cite[Problem~1]{Gorsh_Mas}}]\label{Non-simple_Prob1} Let $G$ be a group which is recognizable by spectrum. What is the largest number $k = k(G)$ such that the direct product of $k$ copies of the group $G$ is still recognizable by spectrum{\rm?}
\end{prob}

\begin{prob}[{\cite[Problem~2]{Gorsh_Mas}}]\label{Non-simple_Prob2} Is it true that for each integer $k \ge 1$ there exists a simple group $G = G(k)$ such that the direct product of $k$ copies of $G$ is recognizable by spectrum{\rm?}
\end{prob}

Problem~\ref{Non-simple_Prob2} was also added by A.~V.~Vasil'ev and V.~D.~Mazurov to Kourovka Notebook~\cite[Problem~20.58a]{Kourovka}, and has been recently solved in the positive in~\cite{YangGorsStarVas}; namely it was proved that if $k$ and $l$ are positive integers such that $n=2^l\ge 56k^2$, then the $k$-th direct power of the group $PSL_n(2)$ is recognizable by spectrum. A survey of some other recent results on recognition by spectrum of groups with non-simple socle and also a number of other open problems in this area can be found in~\cite[Section~4]{GreMazShiVasYang}, see also~\cite[Problem~20.58b]{Kourovka}. In particular, the following theorem holds.

\begin{theorem}[{\cite{Gorshkov,Gorsh_Mas,GreMazShiVasYang,Grech_Vas_abel,Li_Mogh_Vas_Wang,Mazurov97,Wang_ea,YangGorsStarVas}}] \label{DirProd_Spectra} $(i)$ The following non-simple groups are recognizable by spectrum{\rm:}

\begin{itemize}

\item[$(1)$] ${^2}B_2(q)\times {^2}B_2(q)$, where $q=2^l$, where $l \ge 3$ is odd and $l\not = 5$ {\rm (\cite[Theorem~2]{Mazurov97} and \cite[Theorem~3]{Wang_ea})}{\rm;}

\item[$(2)$] the permutation wreath product of ${^2}B_2(2^7)$ with $H \cong 23:11 < S_{23}$ {\rm (\cite[Theorem~1]{Mazurov97})}{\rm;}

\item[$(3)$]  the direct product of the groups $L_i={^2}B_2(2^{p_i})$ for $1\le t\le k$, where $p_1=7$ and for each $i \ge 2$, $p_i$ is the smallest prime lager than $p_{i-1}$ and not lying in $\cup_{j<i}\pi({^2}B_2(2^{p_j}))$ {\rm (\cite[Theorem~4.6]{GreMazShiVasYang})}{\rm;}

\item[$(4)$] $J_4\times J_4$ {\rm (\cite[Theorem]{Gorsh_Mas})}{\rm;}

\item[$(5)$] $J_1\times J_1$  {\rm (\cite[Theorem~1.2]{Li_Mogh_Vas_Wang})}{\rm;}

\item[$(6)$] $PSL_{2^m}(2)\times PSL_{2^m}(2)\times PSL_{2^m}(2)$, where $m>5$  {\rm (\cite[Theorem~1]{Gorshkov})}{\rm;}

\item[$(7)$] the $k$-th direct power of the group $PSL_n(2)$, where $n=2^l\ge 56k^2$  {\rm (\cite[Theorem~1]{YangGorsStarVas})}{\rm;}

\item[$(8)$] ${^2}G_2(q)\times {^2}G_2(q)$ for $q=3^l$, where $l \ge 3$ is odd {\rm (\cite[Theorem~1.1]{Li_Mogh_Vas_Wang})}{\rm;}

\item[$(9)$] the direct product of the groups $L_t={^2}G_2(3^{p_t})$ for $1\le t\le k$, where $p_1=5$ and for each $t \ge 2$, $p_t$ is the smallest prime lager than $p_{t-1}$ and not lying in $\cup_{j<t}\pi({^2}G_2(3^{p_j}))$ {\rm (\cite[Theorem~1]{Grech_Vas_abel})}{\rm;}

\end{itemize}

$(ii)$ The group ${^2}B_2(2^5)\times {^2}B_2(2^5)$ is $4$-recognizable by spectrum {\rm (\cite[Theorem~4]{Wang_ea})}.

$(iii)$ The following non-simple groups are unrecognizable by spectrum{\rm:}

\begin{itemize}

\item[$(1)$]  $PSL_2(q)\times PSL_2(q)$ {\rm (\cite[Proposition~3.4]{Li_Mogh_Vas_Wang})}{\rm;}

\item[$(2)$] $PSp_{2n}(q)\times PSp_{2n}(q)\times PSp_{2n}(q)$, where $q$ is odd {\rm (\cite[Theorem~2]{YangGorsStarVas})}{\rm;}

\item[$(3)$] $PSp_{2n}(q)\times PSp_{2n}(q)$, where $q=p^l$, $p$ is an odd prime, and $2n-1$ is not a power of $p$ {\rm (\cite[Theorem~2]{YangGorsStarVas})}{\rm;}

\item[$(4)$] $PSL_n(q)\times PSL_n(q)\times PSL_n(q)$, where there exists a prime $r$ dividing $q-1$ and coprime to $n$ {\rm (\cite[Theorem~2]{YangGorsStarVas})}{\rm;}

\item[$(5)$] $PSL_n(q)\times PSL_n(q)$,  where $q=p^l$, $p$ is a prime, there exists a prime $r$ dividing $q-1$ and coprime to $n$, and $n-1$ is not a power of $p$ {\rm (\cite[Theorem~2]{YangGorsStarVas})}{\rm;}

\item[$(6)$] $PSU_n(q)\times PSU_n(q)\times PSU_n(q)$, where there exists a prime $r$ dividing $q+1$ and coprime to $n$ {\rm (\cite[Theorem~2]{YangGorsStarVas})}{\rm;}

\item[$(7)$] $PSU_n(q)\times PSU_n(q)$,  where $q=p^l$, $p$ is a prime, there exists a prime $r$ dividing $q+1$ and coprime to $n$, and $n-1$ is not a power of $p$ {\rm (\cite[Theorem~2]{YangGorsStarVas})}{\rm;}

\item[$(8)$] ${^2}G_2(q)\times {^2}G_2(q)\times{^2}G_2(q)$  {\rm (\cite[Theorem~2]{Grech_Vas_abel})}{\rm;}

\item[$(9)$]  $J_1\times J_1\times J_1\times J_1$ {\rm (\cite[Theorem~2]{Grech_Vas_abel})}.

\end{itemize}

\end{theorem}

\section{Characterization by Gruenberg--Kegel graph}

As we discussed in the previous section, ''almost all'' simple groups are almost recognizable by spectrum and there are many non-simple groups which are recognizable by spectrum. Moreover, if $G$ is recognizable by Gruenberg--Kegel graph, then $G$ is recognizable by spectrum; the converse does not hold in general. However, computing the spectrum of a group is not easy. Computing the Gruenberg--Kegel graph of a group is much easier. So, the problem of characterization of a group by its Gruenberg--Kegel graph is of interest.

\smallskip

Recently P.~J.~Cameron and the author~\cite{CaMas} have obtained the following criterion of recognition of a group by Gruenberg--Kegel graph.

\begin{theorem}[{\cite[Theorem~1.2]{CaMas}}]\label{UnrecGKGraph} Let $G$ be a group. The following statements are equivalent{\rm:}
\begin{itemize}
\item[$(1)$] $G$ is unrecognizable by Gruenberg--Kegel graph{\rm;}
\item[$(2)$] there exists a group $H$ such that $S(H)\not =1$ and $\Gamma(H)=\Gamma(G)$.
\end{itemize}
\end{theorem}

It is clear that the implication $(2) \Rightarrow (1)$ follows directly from Theorem~\ref{UnrecSpec} since if $G$ is a group such that there are infinitely many pairwise non-isomorphic groups $H$ with $\omega(G)=\omega(H)$, then each of these groups $H$ satisfies the condition $\Gamma(H)=\Gamma(G)$.

\smallskip

To prove the converse, we use the following scheme. Let $G$ be a group which is unrecognizable by Gruenberg--Kegel graph, and such that if $H$ is a group with $\Gamma(H)=\Gamma(G)$, then $S(H)=1$ (in particular, $H$ is non-solvable). Then $2$ is non-adjacent to at least one odd vertex in $\Gamma(G)$, otherwise $\Gamma(G)=\Gamma(\mathbb{Z}_2\times G)$ and $S(\mathbb{Z}_2\times G)\not = 1$; a contradiction. Therefore by \cite[Proposition~2]{Vasil_2005} (see also Theorem~\ref{vas} below), each $H$ with $\Gamma(H)=\Gamma(G)$ is almost simple, and $\pi(Soc(H))\subseteq \pi(G)$. The following assertion completes the proof.

\begin{prop}[{\cite[Proposition~3.1]{CaMas}}]\label{SimpleNumber} Let $\pi$ be a finite set of primes. The number of non-abelian simple groups $S$ with $\pi(S)\subseteq\pi$ is finite, and is at most $O(|\pi|^3)$.
\end{prop}

Moreover, we have characterized~\cite{CaMas} groups which are almost recognizable by Gruenberg--Kegel graph.

\begin{theorem}[{\cite[Theorem~1.3]{CaMas}}]\label{AlmRecGKGraph} Let $G$ be a group such that $G$ is $k$-recognizable by Gruenberg--Kegel graph for some positive integer $k$. Then the following conditions hold{\rm:}
\begin{itemize}
\item[$(1)$] $G$ is almost simple{\rm;}
\item[$(2)$] each group $H$ with $\Gamma(H)=\Gamma(G)$ is almost simple{\rm;}
\item[$(3)$] each vertex of $\Gamma(G)$ is non-adjacent to at least one other vertex, in particular, $2$ is non-adjacent to at least one odd prime in $\Gamma(G)${\rm;}
\item[$(4)$] $\Gamma(G)$ contains at least $3$ pairwise non-adjacent vertices.
\end{itemize}
\end{theorem}

The next assertion follows directly from Theorem~\ref{AlmRecGKGraph} and Proposition~\ref{SimpleNumber} (see also Theorem~\ref{AlmSimpleNumber} below).

\begin{cor}\label{AlmRecGKGraph_cor} A group $G$ is almost recognizable by Gruenberg--Kegel graph if and only if each group $H$ with $\Gamma(H)=\Gamma(G)$ is almost simple.
\end{cor}

Thus, the problem of recognition of a group by Gruenberg--Kegel graph can be solved in the positive only for almost simple groups since each group which is recognizable by Gruenberg--Kegel graph is almost simple. At the same time, the almost simple groups form a very important class of groups. For example, by the O'Nan--Scott theorem~\cite{ONanScott} (see also~\cite{Asch_Scott}) which was proved in 1988 for primitive permutation groups by M.~W.~Liebeck, C.~E.~Praeger, and J.~Saxl~\cite{LiPraSaxl}, each primitive permutation group is almost simple or is permutation equivalent to a permutation group which belongs to one of described special classes (affine group, simple diagonal action, product action, twisted wreath action)\footnote{Now it is more usual to follow the division and labelling into $8$ classes due to C.~E.~Praeger which appears in~\cite{Praeger,Praeger_Li_Niemeyer}.}.

Moreover, even a very large group can be recognizable by Gruenberg--Kegel graph which can be rather small. For example, the simple sporadic Monster group $M$ with $$|M|=2^{46}\cdot 3^{20}\cdot 5^{9}\cdot 7^{6}\cdot 11^{2}\cdot 13^{3}\cdot 17\cdot 19\cdot 23\cdot 29\cdot 31\cdot 41\cdot 47\cdot 59\cdot 71 \approx 8{,}08\cdot 10^{53}$$ has Gruenberg--Kegel graph pictured in Figure~4.

\begin{center}    \begin{tikzpicture}
        \tikzstyle{every node}=[draw,circle,fill=white,minimum size=4pt, inner sep=0pt]
        \draw (0,0) node (2) [label=below:$2$] {}
        (3.0cm:2.0cm) node (3) [label=above:$3$] {}
        (1.0cm:-1.5cm) node (17) [label=left:$17$] {}
        (0.0cm:-2.0cm) node (7) [label=left:$7$] {}
        (-0.5cm:-2.5cm) node (19) [label=left:$19$] {}
        (-1.5cm:-3.0cm) node (5) [label=above:$5$] {}
        (-2.5cm:-3.0cm) node (11) [label=right:$11$] {}

        (1.0cm:2.5cm) node (31) [label=below:$31$] {}
        (1.5cm:2.5cm) node (23) [label=above:$23$] {}
        (2.0cm:1.5cm) node (13) [label=above:$13$] {}
        (-1.0cm:1.5cm) node (47) [label=above:$47$] {}
        (2.5cm:3.0cm) node (29) [label=below:$29$] {}

        (5.5cm:3.8cm) node (41) [label=above:$41$] {}
        (6.0cm:3.5cm) node (59) [label=above:$59$] {}
        (6.5cm:3.5cm) node (71) [label=above:$71$] {}

        (2)--(17)
        (2)--(7)
        (2)--(19)
        (2)--(5)
        (2)--(11)

        (3)--(17)
        (3)--(7)
        (3)--(19)
        (3)--(5)
        (3)--(11)

        (7)--(17)

        (5)--(11)
        (5)--(19)

        (2)--(3)
        (2)--(47)
        (3)--(29)
        (2)--(13)
        (3)--(13)
        (2)--(23)
        (3)--(23)
        (2)--(31)
        (3)--(31);
    \end{tikzpicture}\\
    Figure~4
    \end{center} By~\cite{Zavarnitsine_2006,LeePopiel}, the group $M$ is recognizable by Gruenberg--Kegel graph. Thus, the following problem naturally arises.

\begin{prob}[{\cite[Problem~1]{CaMas}, see also \cite[Problem~8]{Maslova_conf2023}}]\label{Rec_GK_Prob3} Let $G$ be an almost simple group. Decide whether $G$ is recognizable, $k$-recognizable for some integer $k>1$, or unrecognizable by Gruenberg--Kegel graph.
\end{prob}

Some of the first examples of groups which are recognizable by Gruenberg--Kegel graph were found by M.~Hagie~\cite[Theorem~3]{Hagie_2003} in 2003, namely she proved that sporadic simple groups $J_1$, $M_{22}$, $M_{23}$, $M_{24}$, $Co_2$ are recognizable, and the group $M_{11}$ is $2$-recognizable by Gruenberg--Kegel graph. At the moment there are a number of results on recognition of simple groups by Gruenberg--Kegel graph; a survey of results obtained before $2022$ can be found in~\cite{CaMas} (see also \cite[Remark~1.2, Table~3]{LeePopiel_2} to find updates for recognition question for sporadic simple groups), and in this paper we discuss our recent results in this research area.

In~\cite[Section~7]{Mazurov2005}, V.~D.~Mazurov conjectured that if a simple group $G$ is not isomorphic to $A_6$ and $\Gamma(G)$ has at least three connected components, then $G$ is recognizable by spectrum. This conjecture was proved in a series of papers, the final result was obtained by A.~S.~Kondrat'ev in~\cite{Kondrat2015}, where it was proved that groups $E_7(2)$ and $E_7(3)$ are recognizable by Gruenberg--Kegel graph, so these groups are recognizable by spectrum. In~\cite{MasPansStar} we extended results related to Mazurov's conjecture and considered Problem~\ref{Rec_GK_Prob3} for simple exceptional groups of Lie type whose Gruenberg--Kegel graphs have at least three connected components. The following theorem is true.

\begin{theorem}[{\cite[Main Theorem]{MasPansStar}}]\label{3CompGK} Every simple exceptional group of Lie type, which is isomorphic to neither ${^2}B_2(2^{2n+1})$ with $n\geq1$ nor $G_2(3)$ and whose Gruenberg--Kegel graph has at least three connected components, is almost recognizable by Gruenberg--Kegel graph. Moreover, the groups $ {^2}B_2(2^{2n+1})$, where $n\geq1$, and $G_2(3)$ are unrecognizable by Gruenberg--Kegel graph.
\end{theorem}

In~\cite{Zavarnitsine_2013}, A.~V.~Zavarnitsine proved that if $G=E_8(q)$, where $q \equiv 0, \pm 1 \pmod{5}$ {\rm(}i.\,e., $\Gamma(G)$ has $5$ connected components{\rm)}, and $H$ is a group such that $\Gamma(H)=\Gamma(G)$, then $H \cong E_8(u)$ for a prime power $u \equiv 0, \pm 1\pmod{5}$. In~\cite[Theorem~6.1]{MasPansStar}, we have proved a similar result for groups $G=E_8(q)$, where $q\equiv \pm 2 \pmod{5}$ {\rm(}i.\,e., $\Gamma(G)$ has $4$ connected components{\rm)}. Taking into account Proposition~{\rm\ref{SimpleNumber}} we have immediately that each group $E_8(q)$ is almost recognizable by Gruenberg--Kegel graph. But the question of coincidence of Gruenberg--Kegel graphs of groups $E_8(q)$ and $E_8(q_1)$, where $q\not = q_1$, is still open.

We now turn to the characterization of a group by its order and Gruenberg--Kegel graph. It is known (see, for example, \cite[Propositions~3.1 and~4.3]{VaVd05} and \cite[Proposition~2.4]{VaVd11}) that if $q$ is odd and $n\ge 3$, then $\Gamma(B_n(q))=\Gamma(C_n(q))$ and $|B_n(q)|=|C_n(q)|$ but these groups are not isomorphic. Thus, an analogue of Theorem~\ref{OrdSpecThm} does not hold for the case of characterization by order and Gruenberg--Kegel graph. The following problem was formulated by Behrooz Khosravi in his survey paper~\cite{BKhosravi_survey}, by A.S. Kondrat'ev at the open problems session of the $13$th School–Conference on Group Theory Dedicated to V. A. Belonogov’s 85th Birthday~\cite{Maslova_Conference}, and was independently formulated by W.~Shi in a partial communication with the author.

\begin{prob}[{see \cite[Question~4.2]{BKhosravi_survey}, \cite[Problem~4]{Maslova_Conference}, and \cite[Problem~6.19]{Shi_Surv2024}}]\label{123} For which simple groups $S$ is the following true:
if $G$ is a group with $\Gamma (G)$ = $\Gamma (S)$ and $|G| = |S|$, then $G$
is isomorphic to $S$\/{\rm?}
\end{prob}

It is clear that if a simple group is quasirecognizable by Gruenberg--Kegel graph, then it is uniquely determined by its order and Gruenberg--Kegel graph. The converse does not necessarily hold; for example, in~\cite[Corollary~6.4]{MasPansStar} we have proved that each group $E_8(q)$ is uniquely determined by its order and Gruenberg--Kegel graph while the question of quasirecognition of these groups by Gruenberg--Kegel graph is open.

\section{Approaches to solving Problem~\ref{Rec_GK_Prob3}}

The first step to solve Problem~\ref{Rec_GK_Prob3} is solving the following problem.

\begin{prob}\label{AlmSimpleForRec} Describe the almost simple groups $G$ such that the following statements hold{\rm:}
\begin{itemize}
\item[$(1)$] each vertex of $\Gamma(G)$ is non-adjacent to at least one other vertex, in particular, $2$ is non-adjacent to at least one odd prime in $\Gamma(G)${\rm;}
\item[$(2)$] $\Gamma(G)$ contains at least $3$ pairwise non-adjacent vertices.
\end{itemize}
\end{prob}

Assume that $G$ is an almost simple group such that both statements from Problem~\ref{AlmSimpleForRec} are true and $L=Soc(G)$, and let $H$ be a group such that $\Gamma(H)=\Gamma(G)$. For the question of characterization by Gruenberg--Kegel graph, the cases when the graph is connected and when the graph is disconnected are fundamentally different.

\medskip

Assume first that $\Gamma(G)$ is disconnected. Denote the number of connected components of $\Gamma(G)$
by $s(G)$, and the set of connected components of $\Gamma(G)$ by $\{\pi_i(G) \mid 1 \leq i \leq s(G) \}$; since $|G|$ is even by the Feit-Thompson theorem~\cite[Theorem]{FeitThompson},  we always assume that $2 \in \pi_1(G)$.

The following theorem is a helpful tool for investigations in this case.

\begin{theorem}[{\rm Gruenberg--Kegel Theorem, \cite[Theorem~A]{Williams}}]\label{Gruenberg--Kegel theorem} If~$H$ is a group with disconnected Gruenberg--Kegel graph, then one of the following statements holds{\rm:}
\begin{itemize}
\item[$(1)$] $H$ is a Frobenius group{\rm;}
\item[$(2)$] $H$ is a $2$-Frobenius group{\rm;}
\item[$(3)$] $H$ is an extension of a nilpotent group $N$ by a group~$A$, where $S \unlhd A\le\Aut(S)$,~$S$ is a simple non-abelian group with $s(H)\le s(S)$, $\pi(N)\subseteq \pi_1(H)$, and $\pi(A/S)\subseteq \pi_1(H)$.
\end{itemize}
\end{theorem}

Recall that a group $X$ is called a {\it Frobenius group} if there is a proper non-trivial subgroup $C$ of $X$ such that $C \cap C^x=1$ whenever $x \in X\setminus C$. Let $$
K=\{1_X\} \cup (X \setminus ( \cup_{x \in X} C^x))$$  be the {\em Frobenius kernel} of $X$. It is known~\cite[35.24 and~35.25]{Asch86} that $K \trianglelefteq X$, $X=K \rtimes C$, $C_X(c)\le C$ for each $c \in C$, and $C_X(k)\le K$ for each $k \in K$.
Moreover, by the Thompson theorem on finite groups with fixed-point-free automorphisms of prime order~\cite[Theorem~1]{Thompson}, $K$ is nilpotent. A $2$-Frobenius group is a group $Y$ which contains a normal Frobenius subgroup $X$ with Frobenius kernel $K$ such that $Y/K$ is a Frobenius group with Frobenius kernel $X/K$.

Suppose that $H$ is a Frobenius group or a $2$-Frobenius group. Then $H$ has a non-trivial normal nilpotent subgroup, therefore by Theorem~\ref{UnrecGKGraph}, both $H$ and $G$ are unrecognizable by Gruenberg--Kegel graph. Thus, the following problem naturally arises.

\begin{prob}\label{FrobAlmSimple}
Describe the cases when $\Gamma(G)=\Gamma(H)$, where $G$ is an almost simple group and $H$ is a Frobenius group or a $2$-Frobenius group.
\end{prob}

In the case when $G$ is simple, the complete solution of Problem~\ref{FrobAlmSimple} was obtained by M.~R.~Zinov'eva and V.~D.~Mazurov in~\cite{ZinMaz}. If $H$ is a solvable Frobenius group or a $2$-Frobenius group, then the solution of Problem~\ref{FrobAlmSimple} follows directly from \cite[Theorem~2]{Gorsh_Mas2} and \cite[Theorem~1]{ZinKondr}. The solution of Problem~\ref{FrobAlmSimple} for the case when $H$ is a non-solvable Frobenius group and $Soc(G)\cong PSL_2(q)$ was obtained in the joint note by the author with her PhD student K.~A.~Ilenko~\cite{Masl_Ilenko}, and finally we have considered the remaining cases\footnote{The manuscript is in preparation.}.

Thus, we can suppose that $H$ is an extension of a nilpotent group $N$ by a group~$A$, where $$S \unlhd A\le\Aut(S),$$ $S$ is a simple non-abelian group with $s(H)\le s(S)$, $\pi(N)\subseteq \pi_1(H)$, and $\pi(A/S)\subseteq \pi_1(H)$. Now it is necessary to use the Classification of Finite Simple Groups to go through choices for $S$.

\medskip

In the the case when $\Gamma(G)$ is connected, the situation is more complicated. Here there exists a strong generalization of the Gruenberg--Kegel theorem obtained by A.~Vasil'ev~\cite{Vasil_2005} and then generalized by A.~Vasil'ev and I.~B.~Gorshkov~\cite{VaGo09}.

\medskip

Denote by $t(H)$ the \emph{independence number} of $\Gamma(H)$, that is the greatest size of a coclique  (i.\,e., induced subgraph with no edges) in $\Gamma(H)$. If $r\in\pi(H)$, then denote by $t(r,H)$ the greatest size of a coclique in $\Gamma(H)$ containing $r$.

\begin{theorem}[{\rm \cite[Theorem, Propositions~2 and~3]{Vasil_2005}}]\label{vas}
Let $H$ be a non-solvable group with $t(2,H)\ge 2$. Then the following statements hold.
	
$(1)$ There exists a nonabelian simple group $S$ such that $S \unlhd A = H/N \le\operatorname{Aut}(S)$,  where $N=S(H)$.
	
$(2)$ For every coclique $\rho$ of $\Gamma(H)$ of size at least three, at most one prime in $\rho$ divides the product $|N|\cdot|A/S|$. In particular, $t(S)\geq t(H)-1$.
	
$(3)$ One of the following two conditions holds:	
	
$\mbox{ }$$\mbox{ }$$\mbox{ }$$(3.1)$ $S\cong A_7$ or $L_2(q)$ for some odd $q$, and $t(S)=t(2,S)=3$.
	
$\mbox{ }$$\mbox{ }$$\mbox{ }$$(3.2)$
Every prime $p\in\pi(H)$ nonadjacent to $2$ in $\Gamma(H)$ does not divide the product $|N|\cdot|A/S|$. In particular, $t(2,S)\geq t(2,H)$.
\end{theorem}

Next we discuss the case when $G=L$ is simple. Depending on whether $\Gamma(L)$ is connected or disconnected, we apply Theorem~\ref{Gruenberg--Kegel theorem} or Theorem~\ref{vas}. Further a usual way to prove that a simple group $L$ is recognizable by Gruenberg--Kegel graph is using the following three-step scheme{\rm:}

\begin{itemize}

\item[$(Q)$] prove that $S\cong L$ (i.\,e., prove that $L$ is quasirecognizable by Gruenberg--Kegel graph){\rm;}

\item[$(C)$] prove that $N=1$ (i.\,e., prove that $L$ is recognizable by Gruenberg--Kegel graph among covers){\rm;}

\item[$(A)$] prove that $A/S=1$ (i.\,e., prove that $L$ is recognizable by Gruenberg--Kegel graph among automorphic extensions).

\end{itemize}

Note that to prove that $L$ is almost recognizable by Gruenberg--Kegel graph it is not necessary to prove that $L$ is quasirecognizable  by Gruenberg--Kegel graph (see the paragraph after Theorem~\ref{3CompGK}) or that $L$ is recognizable  by Gruenberg--Kegel among automorphic extensions (since the number of almost simple groups with socle isomorphic to the fixed simple group is finite).

\medskip

\section{On the number of pairwise non-isomorphic groups with the same Gruenberg--Kegel graph}

The question of coincidence of Gruenberg--Kegel graphs of non-isomorphic groups has been studied from different points of view. For example, the author~\cite{Maslova} and independently T.~Burness and E.~Covato~\cite{BurnCov} described all the cases of coincidence of the Gruenberg--Kegel graphs of a simple group and a proper subgroup.

\smallskip

Here we concentrate on the following problem.

\begin{prob} Estimate an upper bound for the number of groups with the same Gruenberg--Kegel graph {\rm(}in the case when this number is finite{\rm)}.
\end{prob}

Let $\Gamma$ be a simple graph whose vertices are labeled by pairwise distinct primes. We call $\Gamma$ a {\it labeled graph}.

In~\cite{CaMas} we proved the following theorem.

\begin{theorem}[{\rm \cite[Theorem~1.4]{CaMas}}]\label{Bounds_GK}  There exists a polynomial function $F(x)=O(x^7)$ such that for each labeled graph $\Gamma$ the following conditions are equivalent{\rm:}
\begin{itemize}
\item[$(1)$] there exist infinitely many pairwise non-isomorphic groups $H$ such that $\Gamma(H)=\Gamma${\rm;}
\item[$(2)$] there exist more than $F(|V(\Gamma)|)$ pairwise non-isomorphic groups $H$ such that $\Gamma(H)=\Gamma$, where $V(\Gamma)$ is the set of the vertices of $\Gamma$.
\end{itemize}

\end{theorem}

Since by Theorem~\ref{AlmRecGKGraph}, if a group $G$ is almost recognizable by Gruenberg--Kegel graph, then each group $H$ with $\Gamma(H)=\Gamma(G)$ is almost simple, a key tool in the proof of Theorem~\ref{Bounds_GK} is the following assertion of independent interest.

\begin{theorem}[{\cite[Theorem~4.2]{CaMas}}]\label{AlmSimpleNumber} There exists a polynomial function $F$ such that if $\pi$ is a finite set of primes, then there are at most $F(|\pi|)$ pairwise non-isomorphic almost simple groups $G$ such that $\pi(G)\subseteq \pi$, and this number is at most $O(|\pi|^7)$.
\end{theorem}

M.~A.~Grechkoseeva and A.~V.~Vasil'ev~\cite[Theorem~3]{GrechVas} generalized our ideas and using their new results on Gruenberg--Kegel graphs of groups with a unique non-abelian composition factor have improved the estimate in Theorem~\ref{Bounds_GK} for the function $F$  to $O(|\pi|^5)$, and the author believes that this result can certainly still  be improved. So, the following problem is of interest.

\begin{prob}[{\cite[Problem~3]{CaMas}}]\label{prob2}
In Theorem~{\rm\ref{Bounds_GK}}, find the exact value for the function $F$, or at least a better upper bound.
\end{prob}

As A.~V.~Zavarnitsine~\cite{Zavarnitsine_2003,Zavarnitsine_2006_1} proved, there is no constant $k$ such that for any almost simple group $G$, the number of pairwise non-isomorphic almost simple groups $H$ such that $\Gamma(G)=\Gamma(H)$ is at most $k$. So, in Theorem~\ref{Bounds_GK} the function $F$ cannot be a constant. However, if $G$ is simple, then A.~V.~Vasil'ev has conjectured that there are at most $4$ simple groups $H$ with $\Gamma(G)=\Gamma(H)$; see Problem~16.26 in~\cite{Kourovka}. This conjecture has been investigated in certain cases. If $G$ is a sporadic simple group, then all simple groups $H$ with $\Gamma(H)=\Gamma(G)$ were described by M.~Hagie in \cite[Corollary~2]{Hagie_2003}. If $G$ is an alternating group, then a similar result was obtained by M.~A.~Zvezdina in \cite{Zvezdina}. The cases of coincidence of Gruenberg--Kegel graphs of groups of Lie type over fields of the same characteristic were described by M.~R.~Zinov'eva~\cite{Zin_1}, and for groups of Lie type over fields of different characteristics partial results were obtained by the same author in a series of papers, a survey can be found in~\cite{Zin_2}.

\section{Characterization by isomorphism type of Gruenberg--Kegel graph}

A group $G$ is {\it recognizable by the isomorphism type of its Gruenberg--Kegel graph} if for each group $H$, $\Gamma(G)$ and $\Gamma(H)$ are isomorphic as abstract graphs (i.\,e., as unlabelled graphs) if and only if $H \cong G$. The following problem naturally arises.

\begin{prob}[{\rm \cite[Problem~13]{Maslova_conf}}]\label{RecIsoGK} Which groups are recognizable by the isomorphism types of their Gruenberg--Kegel graphs{\rm?}
\end{prob}

It is easy to see that if a group $G$ is recognizable by the isomorphism type of its Gruenberg--Kegel graph, then $G$ is recognizable by Gruenberg--Kegel graph, therefore by Theorem~\ref{AlmRecGKGraph}, $G$ is almost simple.

\medskip

There are not many results on recognition of groups by the isomorphism types of their Gruenberg--Kegel graphs. The first result in this direction was obtained by A.~V.~Zavarnitsine~\cite[Theorem~B]{Zavarnitsine_2006}, where it was proved that the simple sporadic group $J_4$ is the unique group whose Gruenberg--Kegel graph has exactly $6$ connected components. In~\cite[Theorem~1.5]{CaMas} we proved that the simple exceptional groups of Lie type ${^2}G_2(27)$ and $E_8(2)$ are also recognizable by the isomorphism types of their Gruenberg--Kegel graphs. Moreover, the author with the participation of M.-Zh.~Chen and M.~R.~Zinov'eva has proved\footnote{The manuscript is in preparation.} that the groups ${^2}E_6(2)$ and $E_8(q)$, where $q \in \{3, 4, 5, 7, 8, 9, 17\}$, are recognizable by the isomorphism types of their Gruenberg–Kegel graphs. Recently Problem~\ref{RecIsoGK} has been solved by M.~Lee and T.~Popiel~\cite{LeePopiel_2} for all simple sporadic groups.

\begin{theorem}[{\cite[Theorem~1.1]{LeePopiel_2}}] A sporadic simple group $G$ is recognisable by the isomorphism type
of its Gruenberg--Kegel graph if and only $G$ is one of the following eight groups{\rm:}
$$B, Fi_{23}, Fi'_{24}, J_4, Ly, M, O'N, Th.$$
Moreover, for every sporadic simple group $G$ not appearing in the above list, there are infinitely many pairwise non-isomorphic groups $H$ with $\Gamma(H)$ isomorphic to $\Gamma(G)$ as abstract graphs.
\end{theorem}

\medskip

So, the following problems arise.

\begin{prob}[{\rm \cite[Problem~5]{CaMas}}] Let $G$ be an almost simple group. Decide whether $G$ is recognizable by the isomorphism type of its Gruenberg--Kegel graph.
\end{prob}

\begin{prob}[{\rm \cite[Comments to Problem~13]{Maslova_conf}}] Which simple groups are recognizable by the isomorphism types of their Gruenberg–Kegel graphs{\rm?}
\end{prob}

\begin{prob}[{\rm \cite[Comments to Problem~13]{Maslova_conf}}] Is there an almost simple but not simple group which is recognizable by the isomorphism type of its Gruenberg–Kegel graph{\rm?}
\end{prob}

\section{Combinatorial properties of Gruenberg--Kegel graphs}

A number of recent papers are devoted to investigation of combinatorial properties of Gruenberg--Kegel graphs of finite groups from different points of view. In this section we will discuss investigations of the following problem.

\begin{prob}\label{UnlabGr} Let $\Gamma$ be a {\rm(}unlabelled{\rm)} graph. Is there a group $G$ such that the graphs $\Gamma(G)$ and $\Gamma$ are isomorphic as abstract graphs {\rm(}i.\,e., as unlabelled graphs{\rm)}{\rm?}
\end{prob}

\medskip
In~\cite{GKMK}, Problem~\ref{UnlabGr} was solved for small graphs. Namely it was proved that for each graph $\Gamma$ with at most $5$ vertices which is not a coclique with $5$ vertices there exists a group $H$ such that $\Gamma(H)$ is isomorphic to $\Gamma$ while a coclique with $5$ vertices is isomorphic to the Gruenberg--Kegel graph of a neither group.
\medskip

The following assertion was first proved in 1999 by M. S. Lusido~\cite{Lucido2}; but the assertion follows directly from the earlier result by G.~Higman~\cite[Theorem~1]{Higman} and the Hall theorem~\cite[Theorem 6.4.1]{Gorenstein_Lect}.

\begin{prop}[{\rm\cite[Proposition~1]{Lucido2}}]\label{non-solvable} Let $G$ be a group with ${t(G) \ge 3}$. Then $G$ is non-solvable.
\end{prop}

A criterion of isomorphism of an arbitrary graph to the Gruenberg--Kegel graph of a solvable group was obtained by A.~Gruber, T.~Keller, M.~Lewis, K.~Naughton, and B.~Strasser in~\cite{Gruber_et_al}. The following theorem holds.

\begin{theorem}[{\rm \cite[Theorem~2]{Gruber_et_al}}]\label{TrianFree1} A graph is isomorphic to the Gruenberg--Kegel graph of a solvable group if and only if its complement is $3$-colorable and triangle-free.
\end{theorem}

The following problem which is a partial case of Problem~\ref{UnlabGr} naturally arises.

\begin{prob}[{\rm \cite[Section~5]{Maslova_CEUR}}]\label{TrianGFree} Is there a graph whose complement is triangle-free but is not $3$-colorable but the graph is isomorphic to the Gruenberg–Kegel graph of an appropriate non--solvable group{\rm?} In the other words, is there a non-solvable group $G$ such that $\Gamma(G)$ does not contain $3$-cocliques and is not isomorphic to the Gruenberg--Kegel graph of any finite solvable group{\rm?}
\end{prob}

Problem~\ref{TrianGFree} was also added by the author to Kourovka Notebook~\cite[Problem~19.52]{Kourovka}. Note that by~\cite[Theorem~1]{Gorsh_Mas2} there are no examples of such groups $G$ from Problem~\ref{TrianGFree} among almost simple groups.

\bigskip

A graph $\Gamma$ is called {\it $k$-regular} if each vertex degree of $\Gamma$ is equal to $k$. A {\it strongly regular graph} with parameters $(v, k, \lambda, \mu)$ is a connected $k$-regular graph $\Gamma$ with $v$ vertices such that every two adjacent vertices have $\lambda$ common neighbours and every two non-adjacent vertices have $\mu$ common neighbours for some integers $\lambda \ge 0$ and $\mu \ge 1$.

The following problem was suggested by Jack Koolen in a private communication with the author.

\begin{prob}[{{\rm J. Koolen, 2016}}] What are the strongly regular graphs which are isomorphic to Gruenberg--Kegel graphs of finite groups{\rm?}
\end{prob}

In~\cite[Theorem~2]{ChenGorsMasYang} M.-Zh.~Chen, I.~B.~Gorshkov, the author and N.~Yang have obtained a solution of this problem. The following theorem holds.

\begin{theorem}[{\rm \cite[Theorem~2]{ChenGorsMasYang}}]\label{SRG}
Let $\Gamma$ be a strongly regular graph such that $\Gamma$ is isomorphic to the Gruenberg--Kegel graph $\Gamma(G)$ of a group $G$. Then one of the following statements holds{\rm:}
	
	$(1)$ $\Gamma$ is the complement to a triangle-free strongly regular graph{\rm;}
	
	$(2)$ $\Gamma$ is a complete multipartite graph with all parts of size $2$.
\end{theorem}

Note that the complement to a complete multipartite graph with all parts of size $2$ is $2$-colorable and triangle-free. Thus, by Theorem~\ref{TrianFree1}, each graph from statement~$(2)$ of Theorem~\ref{SRG} is isomorphic to the Gruenberg--Kegel graph of a solvable group. The following problem which is a partial case of Problem~\ref{TrianGFree} naturally arises.

\begin{prob} Let $\Gamma$ be the complement to a triangle-free but not $3$-colorable strongly regular graph. Is $\Gamma$ isomorphic to the Gruenberg--Kegel graph of a {\rm(}non-solvable{\rm)} group{\rm?}
\end{prob}

Key tools in the proof of Theorem~\ref{SRG} are the following two theorems which are of independent interest.

\begin{theorem}[{\rm \cite[Theorem~1]{ChenGorsMasYang}}]\label{NonAdjTo2} Let $G$ be a finite group of even order such that $t(2,G)\ge 2$. Let $\tau$ be the set of vertices of $\Gamma(G)$ which are not adjacent to $2$. Then the following statements hold{\rm:}
	
	$(1)$ If $G$ is non-solvable, then $G$ has a normal series $$1 \unlhd K \unlhd G_0 \unlhd G,$$ where $K=S(G)$, $G_0/K \cong S$ is a finite non-abelian simple group and $G/K$ is almost simple with socle $S$ and  $$\mbox{either }\tau \subseteq \pi(K) \setminus \pi(G/K) \mbox{ or } \tau \subseteq \pi(S)\setminus (\pi(K) \cup \pi(G/G_0)).$$ In particular, $t(2,G)=2$ or $t(2,G)\le t(2,S)$.
	
	$(2)$ $\tau$ is a union of cliques.
\end{theorem}

Statement~$(1)$ of Theorem~\ref{NonAdjTo2} can be proved by following the arguments in \cite{Vasil_2005} and \cite{VaGo09}. However in~\cite{ChenGorsMasYang} we have found another proof of this statement which is free of considering two different cases. Also note that if $p$ is an odd prime, then there exists a finite group $G$ such that the set of the vertices which are not adjacent to $p$ in $\Gamma(G)$, is connected and is not a clique; see \cite[Remark~1]{ChenGorsMasYang}.

\begin{theorem}[{\rm \cite[Theorem~3]{ChenGorsMasYang}}]\label{ComplMult} Let $\Gamma$ be a complete multipartite graph with each part of size at least~$3$. Then $\Gamma$ is not isomorphic to the Gruenberg--Kegel graph of a group.
\end{theorem}

In~\cite{Masl_Pagon}, the author and D.~Pagon have solved Problem~\ref{UnlabGr} for complete bipartite graphs; namely it was proved that a complete bipartite graph $K_{n,m}$ is isomorphic to the Gruenberg--Kegel graph of a group if and only if $m+n \le 6$ and $(n,m)\not = (3,3)$. So, Theorem~\ref{ComplMult} was proved in~\cite{Masl_Pagon} for the family of complete bipartite graphs.

\section{Acknowledgements}

The author is very thankful to Prof. M.~A.~Grechkoseeva, Prof. V.~V.~Kabanov, Prof. A.~S.~Kondrat'ev, Prof. V.~D.~Mazurov, Prof. C.~E.~Praeger, Prof. W.~Shi, Prof. A.~V.~Vasil'ev, and anonymous reviewers for their helpful comments which improved this manuscript.

\end{document}